\documentclass{birkmult}
\usepackage{amsmath}
\usepackage{amssymb}
\usepackage{url}
\usepackage{esint}

\newcommand{\fun}[6]{\ensuremath{\Phi^{#1,#2}_{#3,#4}(#6,#5)}}

\theoremstyle{plain}
\newtheorem{theorem}{Theorem}[section]
\newtheorem{lemma}[theorem]{Lemma}
\newtheorem{proposition}[theorem]{Proposition}

\theoremstyle{definition}
\newtheorem{definition}[theorem]{Definition}

\theoremstyle{remark}
\newtheorem{remark}[theorem]{Remark}

\begin{document}

\title{A note on boundedness of operators in Grand Grand Morrey spaces }

\author[H. Rafeiro]{Humberto Rafeiro}
\address{Instituto Superior T\'ecnico,\\ Departamento de Matem\'atica, \\Centro CEAF, Av. Rovisco Pais, \\1049--001 Lisboa, \textsc{Portugal}}
\email{hrafeiro@math.ist.utl.pt}

\thanks{H.R. acknowledges the financial support by \textsl{Funda\c c\~ao para a Ci\^encia e a Tecnologia} (FCT), \textsf{Grant SFRH/BPD/63085/2009}, Portugal.}


\dedicatory{Dedicated with great pleasure to Stefan Samko on the occasion of his 70th birthday}

\begin{abstract}
In this note we introduce grand grand Morrey spaces, in the spirit of the  grand Lebesgue spaces. We prove a kind of \textit{reduction  lemma} which is applicable to  a variety of operators to reduce their boundedness in grand grand Morrey spaces to the corresponding  boundedness in  Morrey spaces. As a result of this application, we obtain the boundedness of the Hardy-Littlewood maximal operator and Calder\'on-Zygmund operators in the framework of grand grand Morrey spaces. 
\end{abstract}

\maketitle

\section{Introduction}

In 1992 T. Iwaniec and C. Sbordone \cite{iwa_sbor1992},  in their studies related with the integrability properties of the Jacobian in a bounded open set $\Omega,$ introduced a new type of function spaces  $L^{p)}(\Omega),$   called \textit{grand Lebesgue spaces}. A generalized version of them, $L^{p),\theta}(\Omega)$ appeared in  L. Greco, T. Iwaniec and C. Sbordone \cite{greco}.
 Harmonic analysis related to these spaces and  their associate spaces (called \textit{small Lebesgue spaces}), was intensively studied during last years due to various applications, we mention e.g. \cite{capone_fio, fratta_fio, fio200, fio_gupt_jain, fio_kara, fio_rako, kok_2010}.

Recently in \cite{samko_umar} there was  introduced a version of weighted grand Lebesgue spaces adjusted for sets $\Omega \subseteq \mathbb{R}^n$  of infinite measure, where  the integrability of
$|f(x)|^{p-\varepsilon}$ at infinity was controlled by means of a weight, and there grand grand Lebesgue spaces were also considered, together with the study of classical operators of harmonic analysis in such spaces.  Another idea of introducing ``bilateral" grand Lebesgue spaces on sets of infinite measure was suggested in \cite{357ad}, where the structure of such spaces was investigated, not operators; the spaces in \cite{357ad} are two parametrical with respect to the exponent $p$, with the norm involving $\sup_{p_1<p<p_2}.$ \\

Morrey spaces $L^{p,\lambda}$ were introduced in  1938 by C. Morrey \cite{405a}  in relation to regularity problems of solutions to partial differential equations, and provided a  useful tool in the regularity theory of PDE's (for  Morrey spaces we refer to books \cite{187a, kuf}, see also  \cite{rafsamsam} where an overview of various generalizations may be found).

Recently, in the spirit of grand Lebesgue spaces, A. Meskhi \cite{meskhi2009, meskhi} introduced \textit{grand Morrey spaces}  (in \cite{meskhi2009} it was already defined on quasi-metric measure spaces with doubling measure) and obtained results on the boundedness   of the maximal operator, Cald\'eron-Zygmund singular operators and Riesz potentials. Note that  the ``\textsl{grandification} procedure" was applied only to the parameter $p.$ \\

In this paper we make a further step  and apply the ``grandification procedure" to both the parameters, $p$ and $\lambda,$ obtaining \textit{grand grand Morrey spaces} $L^{p),\lambda)}_{\theta,\alpha}(\Omega).$ In this new framework we obtain a reduction boundedness theorem, which reduces the boundedness of operators (not necessarily linear ones) in grand grand Morrey spaces to the corresponding boundedness in classical Morrey spaces.

\subsection*{Notation}
Throughout the text we use the following notation:

\noindent $\Omega$ stands for an open set in $\mathbb{R}^n$,

\noindent $|A|$ denotes the Lebesgue measure of a measurable set $A\subset \Omega$,

\noindent $B(x,r)=\{y\in \mathbb{R}^n: |y-x|<r\},$ 

\noindent $\widetilde B (x,r) = B(x,r)\cap \Omega,$

\noindent $\mathsf{d}:=\textup{diam}\; \Omega,$

\noindent $\fint_B f(x) \;\mathrm{d}x$ denotes the integral average of the function $f$, i.e. $\fint_B f(x)\;\mathrm{d}x:=|B|^{-1}\int_B f(x)\;\mathrm{d}x,$

\noindent $\hookrightarrow$ means continuous embedding.

\section{Preliminaries}

Everywhere in the sequel, $\Omega$ is supposed to be a bounded open set.  

\subsection{Grand Lebesgue spaces} For $1<p<\infty$, $\theta >0$ and $0<\varepsilon<p-1$ the \emph{grand Lebesgue space} is the set of measurable functions for which
\begin{equation}\label{added_}
 \|f\|_{L^{p),\theta}(\Omega)}:=\sup_{0<\varepsilon<p-1} \varepsilon^\frac{\theta}{p-\varepsilon} \|f\|_{L^{p-\varepsilon}(\Omega)}<\infty.
 \end{equation}
In the case $\theta=1,$ we simply denote $L^{p),\theta}(\Omega):=L^{p)}(\Omega).$

When $|\Omega|<\infty$, then for all $0<\varepsilon\leqslant p-1$ we have
\[
L^p(\Omega)\hookrightarrow L^{p)}(\Omega) \hookrightarrow L^{p-\varepsilon}(\Omega).
\]

For more properties of grand Lebesgue spaces, see \cite{kok_2010}.

\subsection{Morrey spaces}
For $1\leqslant p< \infty $ and $0\leqslant \lambda <1$, the usual Morrey space $L^{p,\lambda}(\Omega)$ is introduced as the set of all measurable functions such that
\[
\|f\|_{L^{p,\lambda}(\Omega)}:=\sup_{\stackrel{x \in \Omega}{0<r\leqslant \mathsf{d}}}\left(\frac{1}{\boldsymbol{|}B(x,r)|^{\lambda}}   \int_{\widetilde{B}(x,r)} |f(y)|^{p} \;\mathrm{d}y\right)^\frac{1}{p}<\infty
\]
where $\mathsf{d}:=\textup{diam}\; \Omega.$

\section{Grand grand Morrey spaces and the reduction lemma}
For $\theta>0$, $\alpha\geqslant0$, $1<p<\infty$ and $0\leqslant \lambda <1$,
 we consider the  functional

 \begin{equation}\label{added}
 \Phi^{p,\lambda}_{\theta,\alpha}(f,s):=\sup_{0<\varepsilon<s} \varepsilon^\frac{\theta}{p-\varepsilon} \|f\|_{L^{p-\varepsilon,\lambda -\alpha\varepsilon}(\Omega)},
 \end{equation}
 where $0<s<\min\{p-1, \lambda/\alpha\}$.
\begin{remark}\label{rem:quotient}
We make a convention that the quotient $\lambda/\alpha$  when $\alpha=0$ is always  $\lambda/\alpha:=\infty$ even if $\lambda=0.$
\end{remark}

\begin{definition}[Grand grand Morrey spaces]
Let $1<p<\infty$, $\theta>0$, $\alpha\geqslant 0$ and $0\leqslant \lambda <1$. By
 $L^{p),\lambda)}_{\theta,\alpha}(\Omega)$ we denote the space of measurable functions having  the finite norm
\begin{equation}\label{equ:norm}
\begin{split}
	\|f\|_{L^{p),\lambda)}_{\theta,\alpha}(\Omega)}:
&=
\Phi^{p,\lambda}_{\theta,\alpha}(f,s_{\max}), \ \ \ \ s_{\max}= \min\left\{p-1,\frac{\lambda}{\alpha}\right\}.
\end{split}
\end{equation}
\end{definition}

\begin{remark}\label{rem:quotient_morrey}
 In the case $\alpha=0, \lambda>0$ we recover  the Grand Morrey spaces introduced in  \cite{meskhi}, and when  $\lambda=\alpha=0$,  by the convention in Remark \ref{rem:quotient} we have the   grand Lebesgue spaces introduced in \cite{greco} (and in \cite{iwa_sbor1992} in the case $\theta=1$).
\end{remark}

For fixed $p,\theta, \lambda, \alpha, f$ we have that $ s \mapsto \Phi^{p,\lambda}_{\theta,\alpha}(f,s)$ is a non-decreasing function, but  it is possible to  estimate $\Phi^{p,\lambda}_{\theta,\alpha}(f,s)$ via $\Phi^{p,\lambda}_{\theta,\alpha}(f,\sigma)$ with $\sigma<s$ as follows.

\begin{lemma}\label{lem:dominance} Let $\Omega$ be a bounded open set.
For $0<\sigma<s<\min\left\{p-1, \frac{\lambda}{\alpha}\right\}$ we have that
\begin{equation}\label{equ:dominance}
	\Phi^{p,\lambda}_{\theta,\alpha}(f,s) \leqslant C s^\frac{\theta}{p-s}\sigma^{-\frac{\theta}{p-\sigma}} \Phi^{p,\lambda}_{\theta,\alpha}(f,\sigma),
\end{equation}
where $C$ depends on   $n,$ the parameters $p,\lambda, \theta, \alpha$  and the diameter $\mathsf d,$ but does not depend on $f, s$ and $\sigma$.
\end{lemma}

\begin{proof}
For fixed $\sigma$ and $0<\sigma<s<\min\{p-1, \frac{\lambda}{\alpha}\}$ we have
\begin{equation}\label{equ:break}
\Phi^{p,\lambda}_{\theta,\alpha}(f,s)= \max \Big\{\Phi^{p,\lambda}_{\theta,\alpha}(f,\sigma) , \underbrace{\sup_{\sigma \leqslant \varepsilon < s} \varepsilon^\frac{\theta}{p-\varepsilon} \|f\|_{L^{p-\varepsilon, \lambda-\alpha\varepsilon}(\Omega)}}_{I} \Big\}.
\end{equation}
To estimate
$$	I=\sup\limits_{\sigma \leqslant \varepsilon < s} \varepsilon^\frac{\theta}{p-\varepsilon}
 \sup\limits_{\stackrel{x \in \Omega}{0<r\leqslant \mathsf{d}}} |B(x,r)|^{\frac{\alpha\varepsilon-\lambda}{p-\varepsilon}} \|f\|_{L^{p-\varepsilon}(\widetilde{B}(x,r))},
 $$
note that the function $g(\varepsilon):=\varepsilon^\frac{\theta}{p-\varepsilon}$ is 
increasing in 
$0<\varepsilon<\min\left\{p-1, \frac{\lambda}{\alpha}\right\},$ 
so that
  \[
\begin{split}	I &\leqslant s^\frac{\theta}{p-s} \sup_{\sigma < \varepsilon < \min\{p-1, \frac{\lambda}{\alpha}\}}\sup_{\stackrel{x \in \Omega}{0<r\leqslant \mathsf{d}}} |B(x,r)|^{\frac{1+\alpha\varepsilon-\lambda}{p-\varepsilon} } \left(  \fint_{\widetilde{B}(x,r)} |f(y)|^{p-\sigma} \;\mathrm{d}y \right)^\frac{1}{p-\sigma}\\
	&\leqslant s^\frac{\theta}{p-s} \sup_{\sigma < \varepsilon < \min\{p-1, \frac{\lambda}{\alpha}\}}\sup_{\stackrel{x \in \Omega}{0<r\leqslant \mathsf{d}}} |B(x,r)|^{\frac{1+\alpha\varepsilon-\lambda}{p-\varepsilon} } \left(  \fint_{\widetilde{B}(x,r)} |f(y)|^{p-\sigma} \;\mathrm{d}y \right)^\frac{1}{p-\sigma}\\
	&\leqslant s^\frac{\theta}{p-s} \!\!   \!\!  \!\!  \!\! \!\! \sup_{\sigma < \varepsilon < \min\{p-1, \frac{\lambda}{\alpha}\}}\sup_{\stackrel{x \in \Omega}{0<r\leqslant \mathsf{d}}} |B(x,r)|^{\Delta(\varepsilon)}
\left( \frac{\sigma^\theta \sigma^{-\theta}}{|B(x,r)|^{\lambda-\alpha\sigma}} \int_{\widetilde{B}(x,r)} |f(y)|^{p-\sigma} \;\mathrm{d}y \right)^\frac{1}{p-\sigma}.
\end{split}
\]
where $\Delta(\varepsilon):=\frac{1+\alpha\varepsilon-\lambda}{p-\varepsilon} -  \frac{1+\alpha\sigma-\lambda}{p-\sigma}$.

Observe that
\[
\Delta(\varepsilon)=\frac{1+\alpha\varepsilon-\lambda}{p-\varepsilon} -  \frac{1+\alpha\sigma-\lambda}{p-\sigma}=\frac{(1-\lambda+\alpha p)(\varepsilon-\sigma)}{(p-\sigma)(p-\varepsilon)}\geqslant 0,
\]
and for $0\leqslant \varepsilon  < \min\{p-1,\lambda/\alpha\}$ we have $\frac{1-\lambda}{p} \leqslant \frac{1+\alpha \varepsilon -\lambda}{p-\varepsilon} \leqslant 1$, so that $0\leqslant \Delta(\varepsilon)\leqslant 1.$ Then
$$|B(x,r)|^{\frac{1+\alpha\varepsilon-\lambda}{p-\varepsilon} -  \frac{1+\alpha\sigma-\lambda}{p-\sigma}}\leqslant C \max\{1, \mathsf d^n\}$$
and we obtain
\[
\begin{split}
	I&\leqslant C s^\frac{\theta}{p-s} \cdot   \sup_{\stackrel{x \in \Omega}{0<r\leqslant \mathsf{d}}}\sigma^{-\frac{\theta}{p-\sigma}} \left( \frac{\sigma^\theta}{|B(x,r)|^{\lambda-\alpha\sigma}} \int_{\widetilde{B}(x,r)} |f(y)|^{p-\sigma} \;\mathrm{d}y \right)^\frac{1}{p-\sigma}\\
&\leqslant Cs^\frac{\theta}{p-s}   \cdot \sigma^{-\frac{\theta}{p-\sigma}}  \sup_{0<\varepsilon \leqslant \sigma}  \varepsilon^\frac{\theta}{p-\varepsilon} \|f\|_{L^{p-\varepsilon,\lambda-\alpha\varepsilon}(\Omega)}\\
&= Cs^\frac{\theta}{p-s} \cdot  \sigma^{-\frac{\theta}{p-\sigma}} \cdot \Phi^{p,\lambda}_{\theta,\alpha}(f,\sigma). \qedhere
\end{split}
\]
\end{proof}

From Lemma \ref{lem:dominance} we immediately have
\begin{lemma}
For $0<\sigma<\min\left\{p-1,\frac{\lambda}{\alpha}\right\}$, the norm defined in \eqref{equ:norm} has the following dominant
\begin{equation}\label{equ:dominant}
	\|f\|_{L^{p),\lambda)}_{\theta,\alpha}(\Omega)}\leqslant C \frac{ \fun{p}{\theta}{\lambda}{\alpha}{\sigma}{f}}{\sigma^\frac{\theta}{p-\sigma}},
\end{equation}
where  $C$ depends on   $n, p,\lambda, \theta, \alpha$  and  $\mathsf d,$ but does not depend on $f$ and $\sigma$.
\end{lemma}

\begin{lemma}[Reduction  lemma]\label{main} Let $U$ be an operator (not necessarily linear) bounded in the usual Morrey spaces
$L^{p-\varepsilon,\lambda-\alpha\varepsilon}(\Omega)$:
	\begin{equation}\label{eq:boun_classical}
\|Uf\|_{L^{p-\varepsilon,\lambda-\alpha\varepsilon}(\Omega)} \leqslant C_{p-\varepsilon,\lambda-\alpha\varepsilon} \|f\|_{L^{p-\varepsilon,\lambda-\alpha\varepsilon}(\Omega)}
\end{equation}
for all sufficiently small $\varepsilon\in [0,\sigma]$, where
 $0<\sigma<\min\left\{p-1,\frac{\lambda}{\alpha}\right\}.$ If we have  $\sup_{0<\varepsilon< \sigma} C_{p-\varepsilon,\lambda-\alpha\varepsilon}<\infty, $ then it is also bounded in the grand grand Morrey space
 $L^{p),\lambda)}_{\theta,\alpha}(\Omega)$:
\begin{equation}\label{equ:metatheorem}
\|Uf\|_{L^{p),\lambda)}_{\theta,\alpha}(\Omega)} \leqslant C\|f\|_{L^{p),\lambda)}_{\theta,\alpha}(\Omega)}
\end{equation}
with
\[
 C=\frac{C_0}{\sigma^\frac{\theta}{p-\sigma}} \sup_{0<\varepsilon< \sigma} C_{p-\varepsilon,\lambda-\alpha\varepsilon},
\]
where $C_0$ may depend on $n,p,\lambda,\theta, \alpha$ and  $\mathsf d,$ but does not depend on $\sigma.$
\end{lemma}

\begin{proof}
By \eqref{equ:dominant},  we have
\begin{equation}\label{equ:bound_operator_grandgrandMorrey}
\|Uf\|_{L^{p),\theta}_{\lambda)}(\Omega)} \leqslant  \frac{C}{\sigma^\frac{\theta}{p-\sigma} } \Phi^{p,\lambda}_{\theta,\alpha}(Uf,\sigma).
\end{equation}
The estimation of $\Phi^{p,\lambda}_{\theta,\alpha}(Uf,\sigma)$ by $\|f\|_{L^{p),\lambda)}_{\theta,\alpha}(\Omega)}$ is direct:
\begin{equation}\label{equ:meta_dominant}
\begin{split}
 \Phi^{p,\lambda}_{\theta,\alpha}(Uf,\sigma)&= \sup_{0<\varepsilon \leqslant \sigma} \varepsilon^\frac{\theta}{p-\varepsilon} \|Uf\|_{L^{p-\varepsilon,\lambda-\alpha\varepsilon}(\Omega)}\\
&\leqslant \sup_{0<\varepsilon \leqslant \sigma} \varepsilon^\frac{\theta}{p-\varepsilon} \cdot C_{p-\varepsilon,\lambda-\alpha\varepsilon} \cdot \|f\|_{L^{p-\varepsilon,\lambda-\alpha\varepsilon}(\Omega)}\\
&\leqslant \sup_{0<\varepsilon \leqslant \sigma}  C_{p-\varepsilon,\lambda-\alpha\varepsilon}\cdot \|f\|_{L^{p),\lambda)}_{\theta,\alpha}(\Omega)}
\end{split}
\end{equation}
which completes the proof.
\end{proof}

\section{On boundedness of operators in the grand grand Morrey spaces}

\subsection{Maximal operator in grand grand Morrey spaces}

Let 
\begin{equation}\label{equ:maximal}
Mf(x)=\sup_{0<r<\mathsf{d}}  \fint_{\widetilde{B}(x,r)} |f(y)| \;\mathrm{d}y, \ \ \ \ x\in\Omega
\end{equation}
be the usual centered maximal operator. The Hardy-Littlewood-Wiener theorem regarding the boundedness of the maximal operator in Lebesgue spaces is a well-known result, see e.g.  \cite{duoandikoetxea}. A similar result is valid in the framework of Morrey spaces, namely


\begin{lemma}\label{lem:maximal_Morrey}
Let $1<p<\infty$ and let $0\leqslant \lambda <1$. Then
\begin{equation}\label{equ:maximal_Morrey}
\|Mf\|_{L^{p,\lambda}(\Omega)} \leqslant \left(2^\frac{n\lambda}{p}c_0 (p^\prime)^\frac{1}{p} +1 \right) \|f\|_{L^{p,\lambda}(\Omega)}.
\end{equation}
\end{lemma}

\begin{remark}
The above explicit evaluation of the constant in Lemma \ref{lem:maximal_Morrey} is the one obtained in \cite{meskhi2009} (see also \cite{kok_2010, meskhi}). For another approach with slightly  different evaluation of the constant, see \cite{alm_has_sam_08, chiafrasca}.
\end{remark}

\begin{theorem}\label{theo:maximal_grandgrandMorrey}
Let $1<p<\infty$, $\theta >0$, $\alpha\geqslant 0$ and $0\leqslant \lambda <1$. Then the Hardy-Littlewood maximal operator \eqref{equ:maximal} is bounded in grand grand Morrey spaces $L^{p),\lambda)}_{\theta,\alpha}(\Omega).$
\end{theorem}

\begin{proof}
By the reduction lemma  \ref{main} and \eqref{equ:maximal_Morrey}, we only need to show the finiteness of
\[
\sup_{0<\varepsilon \leqslant \sigma}C_{p-\varepsilon,\lambda-\alpha\varepsilon}=\sup_{0<\varepsilon \leqslant \sigma} \left(2^\frac{n(\lambda-\alpha\varepsilon)}{p-\varepsilon}c_0 \left(\frac{p-\varepsilon}{p-\varepsilon-1}\right)^\frac{1}{p-\varepsilon} +1 \right)
\]
which holds if we choose  $\sigma<p-1$. Note that the use of the reduction lemma in this proof is not necessary when the grand grand space $L^{p),\lambda)}_{\theta,\alpha}(\Omega)$ is considered with $\alpha>\frac{\lambda}{p-1}.$
\end{proof}

\subsection{Singular integral operators in Grand Grand Morrey spaces}

We follow \cite{meskhi} in this section, in particular, making use of the following definition of the  Calder\'on-Zygmund singular operators.
Namely, the  Calder\'on-Zygmund  operator is treated as the integral operator
\[
Tf(x)=\mathrm{p.v.} \int_\Omega K(x,y)f(y)\;\mathrm{d}y
\]
with the kernel  $K:\Omega \times \Omega \backslash \{(x,x): x \in \Omega\} \to \mathbb R$  satisfying the conditions:
\[
|K(x,y)|\leqslant \frac{C}{|x-y|^n}, \quad x,y \in \Omega, \quad x\neq y;
\]
\[|K(x_1,y)-K(x_2,y)|+|K(y,x_1)-K(y, x_2)| \leqslant Cw\left(\frac{|x_2-x_1|}{|x_2-y|}\right) \frac{1}{|x_2-y|^n}\]
for all $x_1$, $x_2$ and $y$ with $|x_2-y|>C|x-x_2|$, where $w$ is a positive non-decreasing function on $(0,\infty)$ which satisfies the doubling condition $w(2t)\leqslant c w(t)$ and the Dini condition $\int_0^1 w(t)/t \;\mathrm{d}t<\infty$. 
In the case where $w$ is a power function, this goes back to Coifman-Meyers version of  singular operators with \textit{standard kernel}.
We also assume that  $Tf$ exists almost everywhere on $\Omega$  in the principal value sense for all $f \in L^{2}(\Omega)$ and that $T$ is bounded in $L^{2}(\Omega).$ 

The boundedness of such Calder\'on-Zygmund operators in Morrey spaces is valid, as can be seen in the following Proposition, proved in \cite{meskhi}

\begin{proposition}\label{prop:boundedness_ron}
Let $1<p<\infty$, $0\leqslant \lambda <1$. Then
\[
\|Tf\|_{L^{p,\lambda}(\Omega)} \leqslant C_{T,p,\lambda} \|f\|_{L^{p,\lambda}(\Omega)}
\]
where 
\begin{equation}\label{equ:constant_Morrey_calderon}
C_{T,p,\lambda}\le  c\left\{
  \begin{array}{ll}
     \frac{p}{p-1}+\frac{p}{2-p} +\frac{p-\lambda+1}{1-\lambda} & \mbox{if } 1<p<2; \\
    p+\frac{p}{p-2} +\frac{p-\lambda+1}{1-\lambda} & \mbox{if } p>2.\\
  \end{array}
\right.
\end{equation}
with $c$ not depending on $p$ and $\lambda$.
\end{proposition}

\begin{theorem}
Let $1<p<\infty$, $\theta >0$ and $0<\lambda<1$. Then the Calder\'on-Zygmund operator $T$ is bounded in grand grand Morrey spaces $L^{p),\lambda)}_{\theta, \alpha}(\Omega)$.
\end{theorem}

\begin{proof}
  Keeping in mind that by the reduction lemma \ref{main} we are interested only in small values of $\varepsilon$, from  \eqref{equ:constant_Morrey_calderon}, we deduce that
\[
C_{T, p-\varepsilon,\lambda-\alpha\varepsilon}\le 
c \left\{
  \begin{array}{ll}
     \frac{p}{p-\varepsilon-1}+\frac{p-\varepsilon}{2-p+\varepsilon} +\frac{p-\varepsilon-\lambda+\alpha\varepsilon+1}{1-\lambda+\alpha\varepsilon} & \mbox{if } p\le 2 \ \mbox{and} \ \ 0<\varepsilon<p-1; \\
\\
    p-\varepsilon+\frac{p-\varepsilon}{p-\varepsilon-2} +\frac{p-\varepsilon-\lambda+\alpha\varepsilon+1}{1-\lambda+\alpha\varepsilon} & \mbox{if } p>2 \ \ \mbox{and}\ \ 0<\varepsilon<p-2.\\
  \end{array}
\right.
\]
so that when applying the reduction lemma, it suffices to take $\sigma<\min\left\{p-1,\frac{\lambda}{\alpha}\right\}$ when $p\le 2$ and 
$\sigma<\min\left\{p-2,\frac{\lambda}{\alpha}\right\}$ when $p> 2$.
\end{proof}



\subsection*{Acknowledgment}
The author wish to thank Stefan Samko for discussions and useful comments on this paper.
\end{document}